%% file: VZ_question.tex
\documentclass[12pt]{amsart}
\usepackage{amsfonts,amssymb,latexsym}
\usepackage{amsthm,cite,amsmath,amstext,amsbsy,bm,mathrsfs,url} 
\usepackage{paralist}
\usepackage{datetime} 
\usepackage{color}
\usepackage[dvipsnames]{xcolor}
\usepackage[colorlinks,linkcolor=Bittersweet,citecolor=Cerulean]{hyperref} 
\usepackage{mathtools}
\usepackage{euscript}
\usepackage[all]{xy}
\usepackage{dsfont}
\usepackage{graphicx} 
\usepackage[a4paper,text={150mm,230mm},centering]{geometry} 
\usepackage{comment}  
\usepackage{cite}  
\usepackage{thmtools}
\usepackage{enumitem}
\usepackage{letltxmacro} 
\usepackage{nameref}
\usepackage{cleveref}
\usepackage{stmaryrd}
\usepackage{calc}
\usepackage{interval}
\usepackage{manfnt}

\usepackage{tikz-cd}

\makeatletter 
\let\reftagform@=\tagform@
\def\tagform@#1{\maketag@@@{(\ignorespaces\textcolor{PineGreen}{#1}\unskip\@@italiccorr)}}
\renewcommand{\eqref}[1]{\textup{\reftagform@{\ref{#1}}}}
\makeatother

\declaretheorem[
name=Theorem,
Refname={Theorem,Theorems},
numberwithin=section]{thm}
\declaretheorem[
name=Proposition,
Refname={Proposition,Propositions},
sibling=thm]{proposition}

\declaretheorem[
name=Definition,
Refname={Definition,Definitions},
sibling=thm]{dfn}

\newtheorem{thmx}{Theorem}
\renewcommand{\thethmx}{\Alph{thmx}}

\theoremstyle{plain}
\newlist{thmlist}{enumerate}{1}
\setlist[thmlist]{wide = 0pt, labelwidth = 2em, labelsep*=0em, itemindent = 0pt, leftmargin = \dimexpr\labelwidth + \labelsep\relax, noitemsep,topsep = 1ex, font=\normalfont, label=(\roman*), ref=\thethm.(\roman{thmlisti})}

\addtotheorempostheadhook[thm]{\crefalias{thmlisti}{thm}}

\addtotheorempostheadhook[proposition]{\crefalias{thmlisti}{proposition}}

\addtotheorempostheadhook[dfn]{\crefalias{thmlisti}{dfn}}

\addtotheorempostheadhook[lem]{\crefalias{thmlisti}{lem}}
\addtotheorempostheadhook[main]{\crefalias{thmlisti}{main}}

\newlist{thmenum}{enumerate}{1} 
\setlist[thmenum]{wide = 0pt, labelwidth = 2em, labelsep*=0em, itemindent = 0pt, leftmargin = \dimexpr\labelwidth + \labelsep\relax, noitemsep,topsep = 1ex, font=\normalfont, label=(\roman*), ref=\thethmx.(\roman{thmenumi})}
\crefalias{thmenumi}{thmx}

\crefname{lem}{Lemma}{Lemmas}
\crefname{thm}{Theorem}{Theorems}
\crefname{proposition}{Proposition}{Propositions}
\crefname{dfn}{Definition}{Definitions}
\crefname{rem}{Remark}{Remarks}
\crefname{cor}{Corollary}{Corollaries}
\crefname{corx}{Corollary}{Corollaries}
\crefname{problem}{Problem}{Problems}
\crefname{thmx}{Theorem}{Theorems}
\crefname{claim}{Claim}{Claims}

\crefname{main}{Main Theorem}{Main Theorems}

\makeatletter
\newcommand*{\rom}[1]{\expandafter\@slowromancap\romannumeral #1@}
\makeatother

\newcommand{\lowerromannumeral}[1]{\romannumeral#1\relax}

\makeatletter     
\newcommand{\crefnames}[3]{%
	\@for\next:=#1\do{%
		\expandafter\crefname\expandafter{\next}{#2}{#3}%
	}%
}
\makeatother

\crefnames{part,chapter,section}{\S}{\S\S}

\makeatletter
\newsavebox{\@brx}
\newcommand{\llangle}[1][]{\savebox{\@brx}{\(\m@th{#1\langle}\)}%
  \mathopen{\copy\@brx\kern-0.5\wd\@brx\usebox{\@brx}}}
\newcommand{\rrangle}[1][]{\savebox{\@brx}{\(\m@th{#1\rangle}\)}%
  \mathclose{\copy\@brx\kern-0.5\wd\@brx\usebox{\@brx}}}
\makeatother

\let\oldsection\section
\renewcommand{\section}{
	\renewcommand{\theequation}{\thesection.\arabic{equation}}
	\oldsection}
\let\oldsubsection\subsection
\renewcommand{\subsection}{
	\renewcommand{\theequation}{\thesubsection.\arabic{equation}}
	\oldsubsection}

\definecolor{plum}{rgb}{0.8,0.2,0.8}

\newtheorem{rem}[thm]{Remark}

\def\oc{\mathscr{O}}  
  
 \def\fc{\mathcal{F}}

\def\lc{\mathcal{L}}
\def\vc{\mathcal{V}}

\def\cb{\mathbb{C}}

\def\as{{a^\star}}

\def\cb{\mathbb{C}}

\def\as{\mathscr{A}}

\def\fs{\mathscr{F}}

\def\ls{\mathscr{L}}

\def\ts{\mathscr{T}}

\def\dbar{\bar{\partial}}

 \def\d{\partial}

\makeatletter 
\let\@wraptoccontribs\wraptoccontribs
\makeatother

\begin{document}

\title[Brody hyperbolicity of moduli spaces]{Pseudo Kobayashi hyperbolicity of base spaces of  families of minimal projective manifolds with maximal variation}

\author{Ya Deng} 
\email{dengya.math@gmail.com}
\urladdr{https://sites.google.com/site/dengyamath}
	\date{\today} 
	\begin{abstract}
		In this paper we prove that every quasi-projective base space $V$ of   smooth family  of minimal projective manifolds with maximal variation is pseudo Kobayashi hyperbolic, \emph{i.e.}    $V$ is Kobayashi hyperbolic modulo a proper subvariety $Z\subsetneq V$.  In particular, $V$ is algebraically degenerate, that is,    every nonconstant  entire curve  $f:\mathbb{C}\to V$ has image $f(\mathbb{C})$ which lies in that proper subvariety $Z\subsetneq V$. As a direct consequence, we prove the Brody hyperbolicity of moduli spaces of minimal projective manifolds, which answers a question by Viehweg-Zuo in 2003. 
\end{abstract}
\maketitle

\input{section-introduction}
\input{section-higgsbundle}

\bibliographystyle{amsalpha}
\bibliography{biblio}

\end{document}

%% file: section-introduction.tex
\section{Introduction}\label{introduction} 
\subsection{Brody hyperbolicity of moduli spaces}
In \cite{VZ03}  Viehweg-Zuo   proved that moduli spaces of canonically polarized complex manifolds are \emph{Brody hyperbolic}, which can be seen as the  \emph{higher dimensional analytic Shafarevich hyperbolicity conjecture}. 
 Very recently,  Popa-Taji-Wu   \cite{PTW18}  extended Viehweg-Zuo's theorem to moduli spaces of minimal projective manifolds of general type. In this paper, we prove a more general result, which answers a question by Viehweg-Zuo \cite[Question 0.2]{VZ03}.
 \begin{thmx}[Brody hyperbolicity of moduli]\label{VZ question}
 	Let $U\to V$ be a smooth family of polarized   minimal projective manifolds over the quasi-projective manifold $V$ with the Hilbert polynomial $h$. Assume that the induced moduli map $V\to P_h$ is quasi-finite, where $P_h$ denotes to be the (quasi-projective) coarse moduli space associated to the moduli functor $\mathscr{P}_h$. Then $V$ is Brody hyperbolic, \emph{i.e.}   there are no nonconstant entire holomorphic curves on $V$.
 \end{thmx} 

 \subsection{Main tools: Viehweg-Zuo Higgs bundles}
In the series of   works \cite{VZ01,VZ02,VZ03},   Viehweg-Zuo  studied families of minimal projective manifolds   of maximal variation, by constructing the so-called \emph{Viehweg-Zuo Higgs bundles} (see \cref{def:VZ}) on the  base spaces. The Viehweg-Zuo Higgs bundles  (\emph{VZ Higgs bundles} for short) turned to be  a   powerful technique towards understanding the moduli spaces. 
For instance, combing the  celebrated work by Campana-P\u{a}un \cite{CP15,CP15b,CP16}, the existence of VZ Higgs bundles on the aforementioned base spaces  proved in \cite{VZ02} implies the \emph{Viehweg hyperbolicity conjecture}: the base spaces of smooth families of minimal projective manifolds with maximal variation are of log general type. These results were   later extended in \cite{PS17} to bases of families whose generic fibers admit  a good minimal model. When the ``generic Torelli property" (see \cref{main} for the definition) holds for the VZ Higgs bundles, it was established in \cite{VZ03,PTW18} that the base space  is \emph{algebraically degenerate}, \emph{i.e.}  all the non-constant entire curves lie on a proper  subvariety of the base. In \cite{VZ03,PTW18},   such generic Torelli theorems are proved via vanishing theorems, which are unclear to us when the Kodaira dimension of fibers is not maximal. In \cref{main} we  prove that the generic Torelli property is indeed an intrinsic feature  of all VZ Higgs bundles (not related to the Kodaira dimension of fibers!) using different approaches.
\subsection{Pseudo Kobayashi hyperbolicity} 
A quasi-projective variety $V$ is called \emph{pseudo Kobayashi hyperbolic} if $V$ is \emph{Kobayashi hyperbolic modulo a proper subvariety $Z\subsetneq V$}, \emph{i.e.} the Kobayashi-Royden infinitesimal pseudo-metric $\kappa_{ {V}}$ of $V$ is positive definite  outside $Z$. When $Z=\varnothing$, this definition reduces to the usual definition of Kobayashi hyperbolicity.  Such a $V$ is in particular algebraically degenerate. In \cite{Den18},  we proved  that quasi-projective base manifolds of effectively parametrized families of minimal projective manifolds of general type are \emph{Kobayashi hyperbolic}, which generalized the previous work by To-Yeung \cite{TY14}.  A crucial ingredient of our proof is the construction of  (possibly degenerate) Finsler metrics on base spaces whose holomorphic sectional curvatures are  bounded above by a negative constant via VZ Higgs bundles satisfying the generic Torelli property.  Hence by \cite[\S 3]{Den18} and \cref{main}, we conclude the following result.
\begin{thmx}[Pseudo Kobayashi hyperbolicity of bases]\label{Deng}
	Let $f:U\to V$ be a smooth family of minimal projective manifolds over the quasi-projective variety $V$, which is of maximal variation. Then $V$ is pseudo-Kobayashi hyperbolic. In particular, $V$ is algebraically degenerate.
\end{thmx} 
As mentioned above, the Viehweg hyperbolicity conjecture proved by Campana-P\u{a}un shows that the base space $V$ in \cref{Deng} is of log general type.  Hence \cref{Deng} is predicted by  a famous conjecture by S. Lang, which stipulates that a quasi-projective variety
 is pseudo Kobayashi hyperbolic if it is of log general type.
 \begin{rem}
 Using the deep  theories of Hodge modules, Popa-Taji-Wu   constructed VZ Higgs bundles over quasi-projective base spaces of smooth families with maximal variation, whose generic fibers admit   good minimal models (see \cite[Proposition 2.7]{PTW18}). Applying their result  instead of \cref{Higgs bundle} and using our generic Torelli theorem for all VZ Higgs bundles in \cref{main},  \cref{Deng} can be extended to the   pseudo Kobayashi hyperbolicity of those base spaces.
 \end{rem}
 
 \medskip
 
\noindent \textbf{Acknowledgments.}    
 I would like to    thank Professors Junyan Cao,  Jian Xiao and Xiaokui Yang for very helpful discussions, and Professors Damian Brotbek and  Jean-Pierre Demailly  for their  constant supports and encouragements.

%% file: section-higgsbundle.tex
\section{Proofs of  Theorems}\label{construction}
In this section, we first recall  the Viehweg-Zuo Higgs bundles first constructed in  \cite{VZ01,VZ02,VZ03}, and later developed in  \cite{PS17,PTW18}  via Hodge modules. Then we collect some properties (see  also \cite[\S 2]{Den18}), and  prove a \emph{Torelli-type} theorem for the VZ Higgs bundles, which is the main result of our paper.   

\begin{thm}[Viehweg-Zuo, Popa-Taji-Wu]\label{Higgs bundle}
Let $U'\to V'$ be a  smooth family of   minimal projective manifolds over the quasi-projective variety $V'$, which is of \emph{maximal variation}. 
 Then  there exists a birational morphism $\nu:V\to V'$ and a projective compactification  $Y\supset V$, and  two logarithmic Higgs bundles $\big(\bigoplus_{q=0}^{n}F^{n-q,q},\bigoplus_{q=0}^{n}\tau_{n-q,q}\big)$, $\big(\bigoplus_{q=0}^{n}E^{n-q,q},\bigoplus_{q=0}^{n}\theta_{n-q,q}\big)$ together with a big and nef    line bundle \(\ls\)   over $Y$  
  satisfying the following properties:
\begin{thmlist}
	\item \label{dia:higgs}There is a diagram 
\begin{align}\label{dia:two Higgs}
	\xymatrixcolsep{4.3pc}\xymatrix{
	\ls^{-1}\otimes	E^{n-q,q}  \ar[r]^-{\mathds{1}\otimes \theta_{n-q,q}}   &		\ls^{-1}\otimes E^{n-q-1,q+1}\otimes \Omega_Y\big(\log (D+S)\big)  \\
 F^{n-q,q}  \ar[u]^{\rho_{n-q,q}}  \ar[r]^-{  \tau_{n-q,q}}      &     F^{n-q-1,q+1}\otimes \Omega_Y(\log D)   \ar[u]_{\rho_{n-q-1,q+1}\otimes \iota}}
\end{align}
	where both $D:=Y\setminus V$ and $D+S$ are simple normal crossing divisors in $Y$, and $\iota:\Omega_Y(\log D)  \to \Omega_Y(\log D+S)  $ is the natural inclusive map. 
	\item\label{item:2}   $\big(\bigoplus_{q=0}^{n}E^{n-q,q},\bigoplus_{q=0}^{n}\theta_{n-q,q}\big)$  is the logarithmic Higgs bundle underlying the Deligne extension with the real part of the eigenvalues of  residues  in \([0,1)\) of a (polarized) variation of Hodge structures  defined over \(V_0:=Y\setminus D\cup S\).		
	\item \label{injection Higgs} There is an injection $\oc_Y\hookrightarrow F^{n,0}$ which is isomorphic over $V_0$. 
	\item \label{augmented base}  $V_0\cap    \mathbf{B}_+(\ls)=\varnothing$.
	\item \label{iterate}  For any \(k=1,\ldots,n\), \eqref{dia:two Higgs} induces  a  map  
\begin{align}\label{iterated Kodaira}
\tau_k:  {\rm Sym}^k \ts_Y(-\log D)\rightarrow \ls^{-1}\otimes E^{n-k,k}.
\end{align}
\end{thmlist}
\end{thm}
In \cite{VZ02}, the two Higgs bundles in \cref{Higgs bundle} are only  constructed on a big open set of $Y$. In order to extend them to the whole $Y$,  we have to apply a technical lemma in \cite[Proposition 4.4]{PTW18}, so that the  discriminant of the zero divisor defined by  a certain hypersurface used for the cyclic construction, is \emph{a priori normal crossing}.  \cref{iterate} is proved in \cite{VZ03} for canonically polarized families, and  in \cite{PTW18} for smooth families admitting good minimal models.  The readers can also refer to \cite[\S 2]{Den18} for further details.
\begin{proof}[Sketch  of proof of \cref{Higgs bundle}]
By 
\cite[Remark 1.22.(\lowerromannumeral{2})]{Den18} (or \cite[Proposition 3.9]{VZ02} if we are allowed to take a birational model of $U'\to V'$),  after passing to a $r$-folded fiber product $U'^r:=U'\times_{V'}U'\times_{V'}\cdots\times_{V'}U'$,  there exists a projective compactification 
$$ \begin{tikzcd}
 U'^r   \ar[d]  \ar[r,hook] &X' \ar[d, "f'"]      \\
V' \ar[r,hook] &Y'
\end{tikzcd}
$$  
so that
$
mK_{X'/Y'}-mf^*\as
$ is globally generated over $f^{-1}(\Omega)$ for $m\gg 0$, where $\Omega\subset V'$	is some Zariski open set, and $\as$ is a sufficiently ample line bundle. By \cite[Proposition 4.4]{PTW18} (see also \cite[Theorem 1.23]{Den18} if we need to control exceptional locus of the birational morphism), there exists a birational morphism
$$ \begin{tikzcd}
X' \ar[d] &X \ar[l, swap, "\stackrel{{\rm bir}}{\sim}"] \ar[d,"f"]      \\
Y'  & {Y}\arrow{l}{\nu}[swap]{\stackrel{{\rm bir}}{\sim}}  
\end{tikzcd}
$$  
so that there exists a hypersurface
\begin{align}\label{eq:cyclic}
H\in |\ell\Omega_{X/Y}^n(\log \Delta)-\ell f^*\nu^*\as+E|
\end{align} satisfying
\begin{itemize}[leftmargin=0.6cm]
\item there exists a reduced divisor $S$ so that $D+S$ is simple normal crossing, and $H\to Y$ is smooth over $V_0:=Y\setminus (D\cup S)$, where $D:=\nu^{-1}(Y'\setminus V')$. 
\item $E$ is some $f$-exceptional divisor with $f(E)\subset {\rm Supp}(D+S)$. 
\end{itemize} 
Here we denote by $\Delta:=f^{-1}(D)$ so that $(X,\Delta)\to (Y,D)$ is a log morphism.   Within this basic setup, let us define the    two Higgs bundles in the theorem following \cite[\S 4]{VZ02}. 
Leaving out a codimension two subvariety of
${Y}$ supported on $D+S$, we  assume that 
\begin{itemize}[leftmargin=0.6cm] 	
	\item  $f$ is flat, and $E$ in \eqref{eq:cyclic}  disappears. 
	\item The divisor $D+S$ is smooth. Moreover,  both $\Delta $ and $\Sigma=f^{-1}S$ are relative normal crossing. 
\end{itemize}  
 Set $\lc:=\Omega_{X/Y}^n(\log \Delta)$, and $\ls:=\nu^*\as$.  
Let $\delta:W \to X$  be a blow-up of $X$ with centers in $\Delta+\Sigma$ such that  $\delta^*(H+\Delta+\Sigma)$ is a normal crossing divisor. One thus obtains a
cyclic covering of $\delta^*H$,  by taking the $\ell$-th root out of
$\delta^*H$. Let $Z$ 
to be
a strong desingularization of this covering.     We denote the compositions by $h:W\to Y$ and $g:Z\to Y$. Write $\Pi:=g^{-1}(S\cup D)$ which is assumed to be normal crossing.    Leaving out more codimension two subvarieties supported $D+S$, we assume that $h$ and $g$ are also flat, and both $\delta^*(H+\Delta+\Sigma)$ and   $\Pi$ are relative normal crossing.     
Then the restrictions  of both
$g$ and $h$ to
$V_0$ are smooth. Set
$$
F^{n-q,q}:= R^qh_*\Big(\delta^*\big(\Omega_{X/Y}^{n-q}(\log \Delta)\big)\otimes \delta^*\lc^{-1}\otimes \oc_W\big(\lfloor \frac{\delta^*H}{\ell} \rfloor\big)\Big) / {\rm torsion}.
$$ 
It was shown in \cite[\S 4]{VZ02} that there exists a natural  edge
morphism
$$
\tau_{n-q,q}: F^{n-q,q}\to F^{n-q-1,q+1}\otimes \Omega_Y(\log D),
$$
which gives rise to  the first Higgs bundle   $\big(\bigoplus_{q=0}^{n}F^{n-q,q},\bigoplus_{q=0}^{n}\tau_{n-q,q}\big)$ defined over a big open set of $Y$ containing $V_0$.

%
Write $Z_0:=Z\setminus \Pi$.  
Then the local system
$ R^n g_*\cb_{\upharpoonright Z_0}$  
extends to a locally
free sheaf $\vc$ 
on 
$Y$ (here $Y$ is projective rather than the big open set!) equipped with the meromorphic connection
$$
\nabla:\vc\to \vc\otimes \Omega_Y\big(\log (D+S)\big),  
$$ 
whose  eigenvalues of the residues   lie  in $[0,1)$ (the so-called \emph{lower canonical extension}). 
By \cite{Sch73,CKS86}, the Hodge filtration of $ R^n g_*\cb_{\upharpoonright Z_0}$   extends to a filtration $\vc:=\fc^0\supset \fc^1\supset \cdots\supset \fc^{n}$ of \emph{subbundles} so that their graded sheaves $E^{p,n-p}:=\fc^p/\fc^{p+1}$ 
are also  {locally free}, and  there exists
$$
\theta_{p,n-p}:E^{p,n-p}\to E^{p-1,n-p+1}\otimes \Omega_{Y}(\log D+S). 
$$  
This defines the second Higgs bundle 
$
\big(\bigoplus_{q=0}^{n}E^{n-q,q},\theta_{n-q,q}\big)
$. 
As observed in \cite{VZ02,VZ03},   $E^{n-q,q}=R^qg_*\Omega^{n-q}_{Z/Y}(\log \Pi)
$ over a big open set of $Y$ by the Steenbrink's theorem,  
which in turn implies  that  \eqref{dia:two Higgs} holds over a (smaller) big open set in which the Higgs bundle  $\big(\bigoplus_{q=0}^{n}F^{n-q,q},\tau_{n-q,q}\big)$ is also defined,  by \cite[Lemma 6.2]{VZ03} (cf. also \cite[Lemma 4.4]{VZ02}).   

Note that all the objects are defined on a big open set of  $Y$ except for $
\big(\bigoplus_{q=0}^{n}E^{n-q,q},\theta_{n-q,q}\big)
$, which are defined on the whole $Y$.  Following \cite[\S 6]{VZ03}, for every $q=0,\ldots, n$, we define $F^{n-q,q}$ to be 
the reflexive hull, and  the morphism  
$\tau_{n-q,q}$ extends naturally. Since each $E^{n-q,q}$ is locally free, and is defined on   $Y$,  then the morphism $\rho_{n-q,q}$  also
extends. This  leads to \cref{dia:higgs,item:2}.  While \cref{injection Higgs} follows from \cite[Lemma 4.4.(\lowerromannumeral{2})]{VZ02}, \cref{augmented base} can be easily seen from the construction. 

\medskip

To prove \cref{iterate},  we have to introduce a sub-Higgs bundle of $\big(\bigoplus_{q=0}^{n}\ls^{-1}\otimes E^{n-q,q},\bigoplus_{q=0}^{n}\mathds{1}\otimes  {\theta}_{n-q,q}\big)$ following \cite[Corollary 4.5]{VZ02} (or \cite{PTW18}). Write $\tilde{\theta}_{n-q,q}:=\mathds{1}\otimes\theta_{n-q,q}$ for short. For each $q=1,\ldots,n$, we define a coherent torsion free sheaf $\tilde{F}^{n-q,q}:=\rho_{n-q,q}(F^{n-q,q})\subset E^{n-q,q}$ . 
By \cite[Lemma 4.4.(\lowerromannumeral{4})]{VZ02}, $\rho_{n,0}$ is an injection, and thus $\tilde{F}^{n,0}\simeq F^{n,0}\supset \oc_Y$.  By \eqref{dia:two Higgs}, one has
$$
\tilde{\theta}_{n-q,q}:\tilde{F}^{n-q,q}\to \tilde{F}^{n-q-1,q+1}\otimes \Omega_{Y}(\log D),
$$
and let us by $\eta_{n-q,q}$ the restriction of  $\tilde{\theta}_{n-q,q}$ to $\tilde{F}^{n-q,q}$.  Then $\big(\bigoplus_{q=0}^{n}\tilde{F}^{n-q,q},\bigoplus_{q=0}^{n}\eta_{n-q,q}\big)$ is a sub-Higgs bundle of $\big(\bigoplus_{q=0}^{n}\ls^{-1}\otimes E^{n-q,q},\bigoplus_{q=0}^{n}\tilde{\theta}_{n-q,q}\big)$. In particular, $\eta_{n-q+1,q-1}\circ\cdots\circ \eta_{n,0}$ factors like
\begin{align}\label{eq:define iterate2}
\tilde{F}^{n,0}\to  \tilde{F}^{n-q,q}\otimes {\rm Sym}^q\Omega_Y(\log D)\subset \tilde{F}^{n-q,q}\otimes {\bigotimes}^q \Omega_Y(\log D).
\end{align}
By \cref{injection Higgs}, this induces a morphism
\begin{align}\label{eq:define iterate}
\oc_Y\to \tilde{F}^{n,0}\to  \tilde{F}^{n-q,q}\otimes {\rm Sym}^q\Omega_Y(\log D)\hookrightarrow \ls^{-1}\otimes E^{n-q,q}\otimes {\rm Sym}^{q}\Omega_Y(\log D),
\end{align}
and equivalently
$$
\tau_q: {\rm Sym}^{q}\ts_Y(-\log D) \to  \ls^{-1}\otimes E^{n-q,q},
$$
which is the desired morphism in \eqref{iterated Kodaira}.
\end{proof} 
\begin{dfn}\label{def:VZ}
	The negatively twisted Higgs bundle  $ \big(\bigoplus_{q=0}^{n}\ls^{-1}\otimes E^{n-q,q},\bigoplus_{q=0}^{n}\mathds{1}\otimes \theta_{n-q,q}\big)$ satisfying   Properties (\lowerromannumeral{1})-(\lowerromannumeral{5}) in  \cref{Higgs bundle} is called \emph{Viehweg-Zuo Higgs bundle}, and  $\ls$ is called the \emph{Viehweg-Zuo (big) sheaf}.
	\end{dfn}
\begin{thmx}[Generic Torelli property]\label{main}
Same assumption as \cref{Higgs bundle}. Then $\tau_1$ defined in \eqref{iterated Kodaira} is generically injective. 
\end{thmx}
Let us stress here that, differently from \cite[Theorem 2.1.(\lowerromannumeral{6})]{Den18}, we cannot give a precise description of the loci where $\tau_1$ is injective since our method in proving \cref{main} relies on the global aspects of the VZ Higgs bundles. Roughly speaking, the bigness of the Viehweg-Zuo sheaf $\ls$ forces $\tau_1$  to be injective at least one point, which is analogous to Demailly's (weak) holomorphic Morse inequality (see \cite[\S 8.2.(a)]{Dem12}).

Before we prove \cref{main}, we will recall some technical preliminaries.  It was initialed in \cite{VZ03}, and later developed in \cite{PTW18} that, one has to take some proper metric $g$ for $\ls$ so that $g^{-1}$ can compensate the  mild (blow-up) singularities of the Hodge metric $h_{\rm hod}$ of $\bigoplus_{q=0}^{n}E^{n-q,q}$.  \begin{proposition}[\!\!\protect{\cite[Lemma 3.1, Corollary 3.4]{PTW18}}]\label{singular metric}
Same notation as \cref{Higgs bundle}. There exists   a singular metric $g$ for $\ls$ which is smooth over $V_0:=Y\setminus D\cup S$ such that
	\begin{thmlist} 
		\item\label{estimate}  over $V_0$,  the curvature form
	$
		\sqrt{-1}\Theta_{g}(\ls)_{\upharpoonright V_0}$ is positive  definite everywhere.
		\item \label{bounded} The singular hermitian metric \(h:=g^{-1}\otimes h_{\rm hod} \) on \(\ls^{-1}\otimes (\bigoplus_{q=0}^{n}E^{n-q,q})\) is locally bounded  on \(Y\), and smooth over $V_0$. Moreover, \(h\) is \emph{degenerate}  on \(D\cup S \). 
	\end{thmlist}
\end{proposition}
Although the last statement of \cref{bounded} is not explicitly stated in \cite{PTW18}, it can be easily seen from the proof of \cite[Corollary 3.4]{PTW18}. 

\begin{proof}[Proof of \cref{main}]
By \cref{injection Higgs}, $\rho_{n,0}$ induces a global section  $s\in H^0(Y, \ls^{-1}\otimes E^{n,0})$, 
which is \emph{generically} non-vanishing over $V_0$. Set 
\begin{align}\label{set}
V_1:=\{y\in V_0 \mid s(y)\neq 0 \}
\end{align}
which is a non-empty Zariski open set of $V_0$.  For the first stage of VZ Higgs bundle $(\ls^{-1}\otimes E^{n,0},h)_{\upharpoonright V_0}$ over $V_0$, let us denote by $\Theta_0$ its curvature form and set $D'$ to be the $(1,0)$-part of its Chern connection. Then by the Griffiths curvature formula of Hodge bundles (see \emph{e.g.} \cite{GT84}), 	  over $V_0$ we   have
\begin{align}\nonumber
 \Theta_0 &=- \Theta_{\ls,g}\otimes \mathds{1}+\mathds{1}\otimes \Theta_{h_{\rm hod}}(E^{n,0})\\\nonumber
 &=- \Theta_{\ls,g}\otimes \mathds{1}-\mathds{1}\otimes   ({\theta}_{n,0}^*\wedge  {\theta}_{n,0})\\\label{eq:Hodge bundle}
 &=- \Theta_{\ls,g}\otimes \mathds{1}- \tilde{\theta}_{n,0}^*\wedge \tilde{\theta}_{n,0},
\end{align}
where we set \(\tilde{\theta}_{n-k,k}:=\mathds{1}\otimes \theta_{n-k,k} : \ls^{-1}\otimes E^{n-k,k}\rightarrow \ls^{-1}\otimes  E^{n-k-1,k+1}\otimes \Omega_{Y}\big(\log (D+S)\big) \), and define $\tilde{\theta}_{n,0}^*$ to be the adjoint of $\tilde{\theta}_{n,0}$ with respect to the metric $h$. Hence over $V_1$ one has
	\begin{align} \nonumber
	-\sqrt{-1}\d \dbar\log |s|_{h}^2&=  \frac{\big\{ \sqrt{-1}\Theta_0(s),s\big\}_{h}}{| s|^2_{h}}+\frac{\sqrt{-1}\{D's,s \}_h\wedge \{s,D's \}_h}{|s|_{h}^4}-\frac{\sqrt{-1}\{D's,D's \}_h}{|s|_{h}^2} \\\label{eq:crucial}
	& \leqslant   \frac{\big\{ \sqrt{-1}\Theta_0(s),s\big\}_{h}}{| s|^2_{h}}
	\end{align}	
	thanks to the Lagrange's inequality 
	$$\sqrt{-1}|s|^2_{h}\cdot \{D's,D's \}_h\geqslant \sqrt{-1}\{D's,s\}_h\wedge \{s,D's \}_h.$$
Putting \eqref{eq:Hodge bundle} to \eqref{eq:crucial}, over $V_1$ one has 
\begin{align} \label{eq:final}
\sqrt{-1}\Theta_{\ls,g}-\sqrt{-1}\d \dbar\log |s|_{h}^2 \leqslant  -\frac{\big\{ \sqrt{-1}\tilde{\theta}_{n,0}^*\wedge \tilde{\theta}_{n,0}(s),s\big\}_{h}}{| s|^2_{h}}= \frac{\sqrt{-1}\big\{ \tilde{\theta}_{n,0}(s),\tilde{\theta}_{n,0}(s)\big\}_{h}}{| s|^2_{h}} 
\end{align}
where $\tilde{\theta}_{n,0}(s)\in H^0\Big(Y,\ls^{-1}\otimes  E^{n-1,1}\otimes \Omega_{Y}\big(\log (D+S)\big)\Big)$. 
	By \cref{bounded}, for any $y\in D\cup S$, one has
	$$\lim\limits_{y'\in V_0,y'\to y}|s|^2_{h}(y')=0.$$
Therefore, it follows from the compactness of $Y$  that there exists   $y_0\in V_0$ so that $|s|^2_{h}(y_0)\geqslant |s|^2_{h}(y)$ for any $y \in V_0$. Hence $|s|^2_{h}(y_0)>0$, and by \eqref{set}, $y_0\in V_1$. Since $|s|^2_{h}$ is smooth over $V_0$, then
 $
	\sqrt{-1}\d \dbar\log |s|_{h}^2(y_0)
	$ 
	is semi-negative. By \cref{estimate}, $	\sqrt{-1}\Theta_{\ls,g}$ is strictly positive at $y_0$. By \eqref{eq:final} and $|s|_h^2(y_0)>0$, we conclude  that $\sqrt{-1}\big\{  \tilde{\theta}_{n,0}(s),\tilde{\theta}_{n,0}(s)\big\}_{h}$ is strictly positive at $y_0$. In particular, for any non-zero $\xi \in \ts_{Y,y_0}$, $\tilde{\theta}_{n,0}(s)(\xi)\neq 0$. For
	$$
	\tau_1:\ts_Y (-\log D  )\to \ls^{-1}\otimes E^{n-1,1}
	$$
  in \eqref{iterated Kodaira}, over $V_0$ it is defined by $\tau_1(\xi):=\tilde{\theta}_{n,0}(s)(\xi)$ by \eqref{dia:two Higgs}, which 
	is thus \emph{injective at $y_0\in V_1$}. Hence $\tau_1$ is \emph{generically injective}. The theorem is thus proved.
\end{proof}
\begin{proof}[Proof of \cref{Deng}]
By \cref{Higgs bundle}, there exists a VZ Higgs bundle $ \big(\bigoplus_{q=0}^{n}\ls^{-1}\otimes E^{n-q,q},\bigoplus_{q=0}^{n}\mathds{1}\otimes \theta_{n-q,q}\big)$ over a compactification $Y'$ of some birational model $\nu:V'\to V$. Write $D:=Y'\setminus V'$.    In \cite[\S 3]{Den18} we construct a (possibly degenerate) Finsler metric $\mathscr{F}$ over $\ts_{V'}$  defined by 
\begin{align}\label{Finsler}
\fs:=(\sum_{k=1}^{n}{\alpha_k}\fs^2_k)^{1/2}, 
\end{align}
where  $\alpha_1,\ldots,\alpha_n\in \mathbb{R}^+$, and  $\fs_k$ is the Finsler metric on $V'$ defined by
\begin{align}\label{eq:Finslers}
\fs_k(\xi):=|\tau_k(\xi^{\otimes k})|_{h}^{\frac{1}{k}}, \quad \forall \xi \in \ts_{V'}, 
\end{align}
with $\tau_k:{\rm Sym}^k\ts_{Y'}(-\log D)\to \ls^{-1}\otimes E^{n-k,k}$  defined in \eqref{iterated Kodaira}, and $h$  the singular hermitian metric for $ \bigoplus_{q=0}^{n}\ls^{-1}\otimes E^{n-q,q} $ defined in \cref{bounded}. By \cref{main}, $\tau_1$ is injective over   a non-empty Zariski  open set $V''$ of $V'$. Hence $\fs_1$ is positive definite over $V_1:= V''\cap V_0$ by \eqref{eq:Finslers}, where $V_0\subset V'$ is the Zariski open set defined in \cref{Higgs bundle}. By \eqref{Finsler} $\fs$ is also positive definite over the non-empty Zariski open set $V_1$. In \cite[Proposition 3.14]{Den18}, we proved that when $\alpha_1,\ldots,\alpha_k\in \mathbb{R}^+$ are properly chosen, the holomorphic sectional curvature of $\fs$ are bounded from above by a negative constant. It then follows from Demailly's Ahlfors-Schwarz lemma \cite[Lemma 3.2]{Dem97} that the Kobayashi-Royden infinitesimal metric $\kappa_{ {V'}}$ of $V'$ is positive definite over $V_1$. By the \emph{bimeromorphic criteria of Kobayashi hyperbolicity} in \cite[Lemma 3.3]{Den18}, we conclude that $\kappa_{ {V}}$ is positive definite over  the non-empty Zariski open set $\nu(V_1)\subset V$. This proves the theorem. 
\end{proof}
\begin{rem}
	It is natural to ask whether the base space $V'$ of an effectively parametrized family of minimal projective manifolds is Kobayashi hyperbolic. When the fibers are of general type, the answer is positive by \cite[Theorem A]{Den18}. By the proof of \cref{Deng}, in order to prove the Kobayashi hyperbolicity of $V'$, one has to assure that for any given point $p\in V$,
\begin{enumerate}[label=($\spadesuit$),leftmargin=0.6cm]
	\item \label{condition} there exists a birational morphism $\nu:V\to V'$ in \cref{Higgs bundle} which is isomorphic at $p$, and $\nu^{-1}(p)\in V_0$.
\end{enumerate}
\begin{enumerate}[label=($\clubsuit$),leftmargin=0.6cm]
\item \label{condition:Finsler}  $\tau_1$ in \eqref{iterated Kodaira} is injective at the point $\nu^{-1}(p)$.
\end{enumerate} 
 When the fibers are of general type,  \ref{condition} relies on the result of positivity of direct images in \cite[Theorem B.(\lowerromannumeral{3})]{Den18}, and  \ref{condition:Finsler} can be shown via the Bogomolov-Sommese vanishing theorem by  \cite{PTW18}. For general cases, both \ref{condition}  and \ref{condition:Finsler} are unclear to us. Indeed, as mentioned above,  our proof of \cref{main} cannot show the precise loci where $\tau_1$ is injective. 
\end{rem}

A standard inductive arguments in \cite{VZ03} can easily show that \cref{Deng} implies \cref{VZ question}.
\begin{proof}[Proof of \cref{VZ question}]
We will proceed by contradiction. Suppose that there exists a non-constant holomorphic map $f:\cb\to V$. By \cref{Deng}, $f$ cannot be Zariski dense. Let $Z:=\overline{f(\cb)}^{\rm Zar}$ be its Zariski closure. Take a desingularization $\pi:Z'\to Z$. Then there exists a lift $f':\cb\to Z'$ so that $\pi\circ f
'=f$, which is also Zariski dense. Since the new family $U\times_VZ'\to Z'$ is still of maximal variation by the quasi-finiteness of the moduli map,  by \cref{Deng} again, $Z'$ must be algebraically degenerate. This is a contradiction.
\end{proof}